\documentclass{amsart}

\usepackage[latin1]{inputenc}
\usepackage{comment}
\usepackage{qtree,array,youngtab}

\usepackage{color, graphicx}
\usepackage{soul}
\usepackage[T1]{fontenc}

\definecolor{darkgreen}{rgb}{0.,0.5,0.}

\numberwithin{equation}{section}
\newtheorem{theorem}{Theorem}[section]
\newtheorem{lemma}[theorem]{Lemma}

\newtheorem{cor}[theorem]{Corollary}

\allowdisplaybreaks

\title{On the largest size of $\boldsymbol{(t,t+1,\ldots,
t+p)}$-core partitions}

\author{Huan Xiong}
\address{I-Math, Universit$\ddot{a}$t
Z$\ddot{u}$rich, Winterthurerstrasse 190, Z$\ddot{u}$rich 8057,
Switzerland} \email{huan.xiong@math.uzh.ch}



\begin{document}
\begin{abstract}In this paper we prove that Amdeberhan's
conjecture on the largest size of $(t,t+1, t+2)$-core partitions is
true. We also show that the number of $(t, t + 1, t + 2)$-core
partitions with the largest size is $1$ or $2$ based on the parity
of $t$. More generally, the largest size of $(t,t+1,\ldots,
t+p)$-core partitions and the number of such partitions with the
largest size are determined. \end{abstract}

 \maketitle

\section{Introduction}

%background:connection with symmetric groups, p-blocks, why t-cores?

 %how to define partitions and hooks,size, t-core partition?
In number theory and combinatorics, a \emph{partition} is a finite
weakly decreasing sequence of positive integers $\lambda =
(\lambda_1, \lambda_2, \ldots, \lambda_r)$. Let $\mid \lambda
\mid=\sum_{1\leq i\leq r}\lambda_i$. The positive integer $\mid
\lambda \mid$ is called the \emph{size} of the partition $\lambda$.
A partition $\lambda$ could be represented by its Young diagram,
which is a collection of boxes arranged in left-justified rows with
$\lambda_i$ boxes in the $i$-th row. For the $(i, j)$-box, we can
associate its \emph{hook length}, denoted by $h(i, j)$, which is the
number of boxes exactly to the right, or exactly below, or the box
itself. For example, the following are the Young diagram and hook
lengths of the partition $(6, 4, 2)$.

\begin{figure}[htbp]
\begin{center}
\Yvcentermath1

\begin{tabular}{c}
$\young(875421,5421,21)$

\end{tabular}

\end{center}
\caption{The Young diagram of the partition $(6,4,2)$ and the hook
lengths of corresponding boxes. }
\end{figure}

Let $t$ be a positive integer. A partition is called a
\emph{$t$-core partition} if none of its hook lengths is divisible
by $t$. For example, we can see that $\lambda=(6,4,2)$ is a $3$-core
partition from Figure \textbf{$1$}. Furthermore, a partition is
called a $(t_1,t_2,\ldots, t_m)$-core partition if it is
simultaneously  a $t_1$-core, a $t_2$-core, $\ldots$, a $t_m$-core
partition.

 A number of methods, from several areas of
mathematics, have been used in the study of $t$-core and
$(t_1,t_2,\ldots, t_m)$-core partitions. Granville and Ono
\cite{gran} proved that for given positive integers $n$ and $t\geq
4$, there always exists a $t$-core partition with size $n$. It was
showed by Anderson \cite{and} that the number of $(t_1,t_2)$-core
partitions is $\frac{1}{t_1+t_2}  \binom{t_1+t_2}{t_1}$ when $t_1$
and $t_2$ are coprime to each other. Recently, a result obtained by
Olsson and Stanton \cite{ols} was that the largest size of
$(t_1,t_2)$-core partitions is $ \frac{({t_1}^2-1)({t_2}^2-1)}{24}$
when $t_1$ and $t_2$ are relatively prime.

But for general $(t_1,t_2,\ldots, t_m)$-core partitions, what we
know is quite few. We prove the following result, which verifies and
generalizes the conjecture of Amdeberhan \cite{tamd} on the largest
size of $(t,t+1, t+2)$-core partitions:

\begin{theorem} \label{main}
Let $t$ and $p$ be positive integers. Suppose that $t=pn+d$, where
$1\leq d\leq p$ and $n\geq 0.$ Then the largest size of $(t, t + 1,
\ldots, t + p)$-core partitions is
\begin{eqnarray*}
&\max\{\binom{n+2}{2}[\frac{d}{2}](d-[\frac{d}{2}])+\binom{n+2}{3}(p^2n+pd-p^2)-3\binom{n+2}{4}p^2,&\\
&\binom{n+1}{2}(p-[\frac{p-d}{2}])(d+[\frac{p-d}{2}])+\binom{n+1}{3}(p^2n+pd-p^2)-3\binom{n+1}{4}p^2
\},&\end{eqnarray*} where $\max\{x,\ y\}$ denotes the maximal
element in $\{x,\ y\}$. The number of $(t, t + 1, \ldots, t +
p)$-core partitions with the largest size is at most $4$.
\end{theorem}

\begin{cor}
(Cf. Conjecture 11.2 of  \cite{tamd}.) The largest size $g(t)$ of
$(t, t + 1, t + 2)$-core partitions equals to:

\[ g(t) =
\begin{cases}
\  \  \  \  \  \ \ n \binom{n+1}{3}, \ \ \ \ \ \ \ \ \ \ \    \text{if} \ t = 2n-1;\\
  ( n+1) \binom{n+1}{3}+\binom{n+2}{3},\ \text{if} \ t = 2n.
\end{cases}
\]
\end{cor}

\section{The $\boldsymbol{\beta}$-sets of $\boldsymbol{(t,t+1,\ldots,
t+p)}$-core partitions}

 Let $\lambda = (\lambda_1, \lambda_2, \ldots, \lambda_r)$ be a partition whose corresponding Young diagram has $r$ rows.
  The \emph{$\beta$-set} of the partition $\lambda$ is
 denoted by
$$\beta(\lambda)=\{h(i,1) : 1 \leq i \leq r\},$$ which is the set of
hook lengths of boxes in the first column of the corresponding Young
diagram. It is easy to see that $h(1,1)> h(2,1)>\cdots>h(r,1)>0$ and
thus $\beta(\lambda)\subseteq \{0,1, 2,  \ldots , h(1,1)\}$.
 Let $\beta(\lambda)'$ be the complement of
$\beta(\lambda)$ in $\{0,1, 2,  \ldots , h(1,1)\}$ and $H(\lambda)$
be the multiset of hook lengths of $\lambda$. Then
$\beta(\lambda)\subseteq H(\lambda).$ We know $0\in \beta(\lambda)'$
since $0\notin \beta(\lambda)$. It is easy to see that $\lambda$ is
a $t$-core partition if and only if $H(\lambda)$ doesn't contain any
multiple of $t$. The following results are well-known and easy to
prove:

\begin{lemma} \label{thm1}
(\cite{berge})  The  partition $\lambda$ is uniquely determined by
its $\beta$-set.

(1) Suppose $\lambda = (\lambda_1, \lambda_2, \ldots, \lambda_r)$.
Then $\lambda_i=h(i,1)-r+i$ for $ 1 \leq i \leq r$. Thus the size of
$\lambda$ equals to $\mid \lambda \mid=\sum_{x\in
\beta(\lambda)}{x}-\binom{ \#\beta(\lambda) }{2}$, where $\#
\beta(\lambda) $ denotes the number of elements in $\beta(\lambda)$;

(2) $H(\lambda) = \{x - x' : x\in \beta(\lambda), \ x' \in
\beta(\lambda)' ,\ x
> x'\}$.
\end{lemma}

\noindent\textbf{Remark.} Any finite set of some positive integers
could be a $\beta$-set of some partition. Actually, by Lemma
\ref{thm1}, it is easy to see that, given any finite set $S$ of some
positive integers, we can recover a partition by considering $S$ as
a $\beta$-set. Then we know there is a bijection between partitions
and finite sets of some positive integers.

Any finite positive integer set could be a $\beta$-set of some
partition. But to be a $\beta$-set of some $t$-core partition, a
finite positive integer set must satisfy the following condition.

\begin{lemma} \label{a-mt}
A partition $\lambda$ is a $t$-core partition if and only if for any
$x\in \beta(\lambda)$ and any positive integer  $m$ with $x\geq mt$,
we have $x-mt \in \beta(\lambda)$.
\end{lemma}
\textbf{Proof.} $\Rightarrow$: Suppose that $\lambda$ is a $t$-core
partition, $x\in \beta(\lambda)$,\ $m$ is a positive integer, and
$x\geq mt$. By the definition of t-core partitions, we have $mt
\notin H(\lambda) $ and thus $x>mt$. But we know $x-(x-mt)=mt \notin
H(\lambda)$, $x\in \beta(\lambda)$, and $x>x-mt$. Then by Lemma
\ref{thm1}(2), $x-mt$ couldn't be an element in $\beta(\lambda)'$.
Thus we know $x-mt \in \beta(\lambda)$.

$\Leftarrow$: Suppose that for any $x\in \beta(\lambda)$ and any
positive integer  $m$ with $x\geq mt$, we have $x-mt \in
\beta(\lambda)$. This means that for any such $x$ and $m$ we have
$x-mt\notin \beta(\lambda)'$. Thus for any $x\in \beta(\lambda), \
x' \in \beta(\lambda)', \ x
> x'$, we
know $x-x'$ couldn't be a multiple of $t$. Then by Lemma
\ref{thm1}(2) we know $\lambda$ must be a $t$-core partition. \hfill
$\square$

Throughout this paper, let $t$ and $p$ be positive integers. We have
the following lemmas.

\begin{lemma} \label{rep}
Let $k$ be a positive integer. Then \begin{eqnarray*} && \{\sum
\limits_{0\leq i \leq p}{c_i(t+i)}: c_i\in \textbf{Z},\ c_i\geq 0\
(0\leq i \leq p),\ \sum \limits_{0\leq i \leq p}{c_i}=k\}\\ &=& \{ \
x\in \textbf{Z}: kt\leq x \leq k(t+p) \}.\end{eqnarray*}
\end{lemma}
\textbf{Proof.} Suppose that $c_i\in \textbf{Z},\ c_i\geq 0\ (0\leq
i \leq p)$ and $\sum_{0\leq i \leq p}{c_i}=k.$ Let $x=\sum_{0\leq i
\leq p}{c_i(t+i)}$. It is easy to see that
$$kt=\sum \limits_{0\leq i \leq p}{c_it}\leq x \leq \sum
\limits_{0\leq i \leq p}{c_i(t+p)}= k(t+p).$$

On the other hand, suppose that $x\in \textbf{Z}$ and $kt\leq x \leq
k(t+p)$. We will show by induction that $$x\in \{\sum \limits_{0\leq
i \leq p}c_i(t+i): c_i\in \textbf{Z},\ c_i\geq 0\ (0\leq i \leq p),\
\sum \limits_{0\leq i \leq p}c_i=k\}.$$ First it is obvious that
$$kt \in \{\sum \limits_{0\leq i \leq p}c_i(t+i): c_i\in
\textbf{Z},\ c_i\geq 0\ (0\leq i \leq p),\ \sum \limits_{0\leq i
\leq p}c_i=k\}.$$ Suppose that for  $kt< x\leq k(t+p) $ we already
have $x-1=\sum_{0\leq i \leq p}c_i(t+i)$ for some $c_i\in
\textbf{Z},\ c_i\geq 0\ (0\leq i \leq p)$  and $ \sum_{0\leq i \leq
p}c_i=k$. Now we have $c_p< k$ since $x-1< k(t+p)$. Then there must
exist some $0\leq i_0 \leq p-1$ such that $c_{i_0}\geq 1$. Thus we
have \begin{eqnarray*} x&=& 1+ \sum \limits_{0\leq i \leq
p}c_i(t+i)\\ &=&\sum \limits_{0\leq i \leq p,\ i\neq i_0,\
i_0+1}c_i(t+i)+(c_{i_0}-1)(t+i_0)+(c_{i_0+1}+1)(t+i_0+1).\end{eqnarray*}
 It follows
that $$x\in \{\sum \limits_{0\leq i \leq p}c_i(t+i): c_i\in
\textbf{Z},\ c_i\geq 0\ (0\leq i \leq p),\ \sum \limits_{0\leq i
\leq p}c_i=k\}.$$ Now we finish the induction and  prove the lemma.
\hfill $\square$

\begin{lemma} \label{linear}
Let $\lambda$ be a $(t,t+1,\ldots, t+p)$-core partition. Suppose
that $c_i\in \textbf{Z}$  and  $ c_i\geq 0$ for $ 0\leq i \leq p$.
Then $\sum_{0\leq i \leq p}c_i(t+i) \notin \beta(\lambda)$.
\end{lemma}
\textbf{Proof.} Let $k=\sum_{0\leq i \leq p}c_i$. We will prove this
lemma by induction on $k$.  If $k=0$, we have $\sum_{0\leq i \leq
p}c_i(t+i)=0 \notin \beta(\lambda)$. Now assume that $k\geq 1$ and
the result is true for $k-1$. Assume that the result is not true for
$k$, i.e.,  there exist $c_i\in \textbf{Z},\ c_i\geq 0\ (0\leq i
\leq p)$ such that $\sum_{0\leq i \leq p}c_i=k$ and $\sum_{0\leq i
\leq p}c_i(t+i) \in \beta(\lambda)$. Then there must exist some
$0\leq i_0 \leq p$ such that $c_{i_0}\geq 1$ since $\sum_{0\leq i
\leq p}c_i=k\geq 1$. By Lemma \ref{a-mt}, $$\sum \limits_{0\leq i
\leq p}c_i(t+i)-(t+i_0)= \sum \limits_{0\leq i \leq p,\ i\neq
i_0}c_i(t+i)+(c_{i_0}-1)(t+i_0)\in \beta(\lambda)$$ since $\lambda$
is a $(t,t+1,\ldots, t+p)$-core partition. But by assumption we know
$$\sum_{0\leq i \leq p,\ i\neq i_0}c_i(t+i)+(c_{i_0}-1)(t+i_0)\notin
\beta(\lambda)$$ since $\sum_{0\leq i \leq p,\ i\neq
i_0}c_i+(c_{i_0}-1)=k-1$, a contradiction! This means that we must
have $$\sum_{0\leq i \leq p}c_i(t+i) \notin \beta(\lambda)$$ for
$c_i\in \textbf{Z},\ c_i\geq 0\ (0\leq i \leq p)$ and $\sum_{0\leq i
\leq p}c_i=k$. We finish the induction. \hfill $\square$

Let $[ x ]$ be the largest integer not greater than $x$. For $1\leq
k \leq [\frac{t+p-2}{p}]$, let $$S_{k}=\{  x\in \textbf{Z} :
(k-1)(t+p)+1\leq x \leq  kt-1 \}.$$

Notice that  for $1\leq k \leq [\frac{t+p-2}{p}]$,
$S_k\neq\emptyset$ since $(k-1)(t+p)+1\leq  kt-1$.

We have the following characterization for $\beta$-sets of
$(t,t+1,\ldots, t+p)$-core partitions.

\begin{lemma} \label{set}
Let $t$ and $p$ be positive integers. Suppose that $\lambda$ is a
$(t,t+1,\ldots, t+p)$-core partition.  Then $\beta(\lambda)$ must be
a subset of $ \  \bigcup_{1\leq k \leq [\frac{t+p-2}{p}]} S_{k}$.
\end{lemma}
\textbf{Proof.} First we claim that for every $x\geq
[\frac{t+p-2}{p}]t$, we have $x\notin \beta(\lambda)$:

Suppose $x\geq [\frac{t+p-2}{p}]t$. Then there must exist some
$k\geq [\frac{t+p-2}{p}]$ such that $kt\leq x < (k+1)t$. Thus we
know  $$ kt\leq x \leq  (k+1)t-1\leq  kt+   [\frac{t+p-2}{p}]p\leq
k(t+p)$$ since $k\geq [\frac{t+p-2}{p}]$. By Lemma \ref{rep} we have
$$x\in \{\sum \limits_{0\leq i \leq p}c_i(t+i): c_i\in \textbf{Z},\
c_i\geq 0\ (0\leq i \leq p),\ \sum \limits_{0\leq i \leq
p}c_i=k\}.$$ Then by Lemma \ref{linear} we know $x\notin
\beta(\lambda)$. The claim is proved.

Now we know $\beta(\lambda)$ must be a subset of $\{x\in \textbf{Z}
: 1\leq x\leq [\frac{t+p-2}{p}]t-1 \}$. By Lemma \ref{rep} and Lemma
\ref{linear} we have $$\{ \ x\in \textbf{Z}: kt\leq x \leq k(t+p) \}
\bigcap \beta(\lambda)=\emptyset$$ for every positive integer $k$.
Hence $\beta(\lambda)$ must be a subset of $$\{
 x\in \textbf{Z} :1\leq x\leq [\frac{t+p-2}{p}]t-1 \} \setminus (\mathop{\bigcup}_{1\leq k
\leq [\frac{t+p-2}{p}]-1} {\{x\in \textbf{Z}: kt\leq x \leq k(t+p)
\}}),$$ which equals to $\bigcup_{1\leq k \leq [\frac{t+p-2}{p}]}
S_{k}$. \hfill $ \square$

By Lemma \ref{set} and Lemma \ref{thm1}, the next result is obvious.
We mention that, the following result is also a corollary of Theorem
\textbf{$1$} in \cite{and}.
\begin{cor} \label{}
Let $t$ and $p$ be positive integers. Then the number of
$(t,t+1,\ldots, t+p)$-core partitions must be finite.
\end{cor}

\section{The largest size of $\boldsymbol{(t, t+1, \ldots, t+p)}$-core partitions}
%In this section , we will prove the largest size of $(t,t+1,\ldots,
%t+p)$-core partition is $???? $ and the number of
% $(t,t+1,\ldots, t+p)$-core partition who has this largest
%size is $t-2[\frac{t}{2}]+1$.

Let $\lambda$ be a $(t,t+1,\ldots, t+p)$-core partition. By Lemma
\ref{set}, we know $\beta(\lambda)\subseteq \bigcup_{1\leq k \leq
[\frac{t+p-2}{p}]} S_{k}$. Let $a_{k}=\# S_{k} $ be the number of
elements in $S_k$ and $b_{\lambda,k}=\# (\beta(\lambda)\bigcap
S_{k})
  $ be the number of elements in $\beta(\lambda)\bigcap S_{k}
  $ for $1\leq k \leq [\frac{t+p-2}{p}]$. It is obvious that $b_{\lambda,k}\leq
  a_k$ for $1\leq k \leq [\frac{t+p-2}{p}]$.

\begin{lemma} \label{a_k}
Let $1\leq k \leq [\frac{t+p-2}{p}] $. Then $a_{k}=t-(k-1)p-1$. Thus
 $1\leq a_{[\frac{t+p-2}{p}]}\leq p$ and  for $1\leq k \leq
[\frac{t+p-2}{p}]-1 $,  we have $a_{k}-a_{k+1}=p$. Additionally, for
every $x\in S_{k+1}$ and $0\leq i\leq p$, we have $x-(t+i)\in
S_{k}$. Furthermore, $\bigcup_{1\leq k \leq [\frac{t+p-2}{p}]}
S_{k}$ is a $\beta$-set of some $(t,t+1,\ldots, t+p)$-core
partition.
\end{lemma}
\textbf{Proof.} First we have
$$a_{k}=kt-1-((k-1)(t+p)+1)+1=t-(k-1)p-1$$ for $1\leq k \leq
[\frac{t+p-2}{p}] .$ Thus
$$1\leq t-([\frac{t+p-2}{p}]-1)p-1=a_{[\frac{t+p-2}{p}]}\leq
t-(\frac{t+p-2-(p-1)}{p}-1)p-1= p$$ and
$$a_{k}-a_{k+1}=t-(k-1)p-1-(t-kp-1)=p.$$

 Suppose that $x\in S_{k+1}$
for $1\leq k \leq [\frac{t+p-2}{p}]-1 $. This means that
$$k(t+p)+1\leq x\leq (k+1)t-1.$$ Thus for $0\leq i\leq p$, we have
$$(k-1)(t+p)+1\leq k(t+p)+1-(t+i)\leq x-(t+i)\leq (k+1)t-1-(t+i)\leq
kt-1,$$ which means that $x-(t+i)\in S_{k}$. Then by Lemma
\ref{a-mt}  we know $\bigcup_{1\leq k \leq [\frac{t+p-2}{p}]} S_{k}$
must be a $\beta$-set of some $(t,t+1,\ldots, t+p)$-core partition.
\hfill $\square$

\begin{lemma} \label{bk}
 Let $\lambda$ be a $(t,t+1,\ldots, t+p)$-core partition and $1\leq k \leq
[\frac{t+p-2}{p}]-1 $. If $b_{\lambda,k+1}\neq 0$, then
 $b_{\lambda,k}-b_{\lambda,k+1}\geq p$.
\end{lemma}
\textbf{Proof.} Suppose that $1\leq k \leq [\frac{t+p-2}{p}]-1 $ and
$b_{\lambda,k+1}\neq 0$. Let
$$\beta(\lambda)\bigcap S_{k+1}=\{x_i: 1\leq i \leq
b_{\lambda,k+1}\}$$ where $x_1 < x_2< \cdots <x_{b_{\lambda,k+1}}$.
Then by Lemma \ref{a-mt}
 and Lemma \ref{a_k} we know $\{ x_1-(t+p), x_1-(t+p-1),\ldots, x_{1}-(t+1),  x_1-t, x_2-t,
\ldots,  x_{b_{\lambda,k+1}}-t \}\subseteq \beta(\lambda)\bigcap
S_{k}$ where $x_1-(t+p)< x_1-(t+p-1)<\cdots< x_{1}-(t+1)< x_1-t<
x_2-t< \cdots< x_{b_{\lambda,k+1}}-t$. It follows that
$\beta(\lambda)\bigcap S_{k}$ has at least $b_{\lambda,k+1}+ p$
different elements and thus $b_{\lambda,k}\geq b_{\lambda,k+1}+ p$.
 \hfill $\square$

Let $1\leq r \leq [\frac{t+p-2}{p}].$ Suppose that $c_1, c_2,
\ldots, c_r$ are positive integers and $c_k\leq a_{k}$ for $1\leq k
\leq r$. Let $\mu_{c_1, c_2, \ldots, c_r}$ be the partition whose
$\beta$-set satisfies $$\beta(\mu_{c_1, c_2, \ldots, c_r})\subseteq
\bigcup_{1\leq k \leq r} S_{k}$$    and
$$\beta(\mu_{c_1, c_2, \ldots, c_r})\bigcap S_{k}= \{ x\in
\textbf{Z} : kt-c_k\leq x \leq kt-1 \}$$ for $1\leq k \leq r.$

\begin{lemma} \label{mostnumbers}
Suppose that $c_k\leq a_{k}$ for $1\leq k \leq r$. The partition
$\mu_{c_1, c_2, \ldots, c_r}$ is a $(t,t+1,\ldots, t+p)$-core
partition if and only if $c_k-c_{k+1}\geq p$ for $1\leq k \leq r-1$.
\end{lemma}
\textbf{Proof.}    Suppose that $1\leq k \leq r-1$ and $x\in
\beta(\mu_{c_1, c_2, \ldots, c_r})\bigcap S_{k+1}.$ This means that
$(k+1)t-c_{k+1}\leq x \leq (k+1)t-1$. Thus for $0\leq i\leq p$, we
have
$$ (k+1)t-c_{k+1}-(t+p)\leq x-(t+i) \leq (k+1)t-1-t=kt-1.$$
Then by Lemma \ref{a-mt} and Lemma \ref{a_k} it is easy to see that
$\mu_{c_1, c_2, \ldots, c_r}$ is a $(t,t+1,\ldots, t+p)$-core
partition if and only if $kt-c_k\leq (k+1)t-c_{k+1}-(t+p)$ for
$1\leq k \leq r-1$, which is equivalent to $c_k-c_{k+1}\geq p$ for
$1\leq k \leq r-1$. \hfill $\square$

Let ${\gamma}_i= \mu_{i, i-p, i-2p,  \ldots, i- [\frac{i-1}{p}] p}$
and $f(i)=\mid {\gamma}_i  \mid$ be the size of ${\gamma}_i$ for
$1\leq i \leq t-1$.   By Lemma \ref{mostnumbers} ${\gamma}_i$ is a
$(t,t+1,\ldots, t+p)$-core partition for $1\leq i \leq t-1$. For
convenience, let ${\gamma}_0$ be the empty partition and $f(0)=0$.

\begin{lemma} \label{f(pm+i+1)-f(pm+i)}
Suppose that $1\leq i \leq p$, $m\geq 0$ and $pm+i\leq t-1.$ Then we
have

(1)  $f(pm+i)-f(pm+i-1)=\binom{m+2}{2}(t-pm-2i+1)$;

(2) $ f(pm+i)-f(pm)=\binom{m+2}{2}(it-ipm-i^2). $
\end{lemma}
\textbf{Proof.} \textbf{(1)} First we know
${\gamma}_{pm+i}=\mu_{pm+i, p(m-1)+i, p(m-2)+i,  \ldots, i }$ and
$$ \#
\beta({\gamma}_{pm+i})= \sum \limits_{0\leq j \leq
m}(pj+i)=p\binom{m+1}{2}+(m+1)i.$$ Then by Lemma \ref{thm1}(1) we
have
\begin{eqnarray*}
& &f(pm+i)-f(pm+i-1)\\ &=& \mid   {\gamma}_{pm+i}  \mid- \mid
{\gamma}_{pm+i-1} \mid  \\ &=& \sum \limits_{x\in
\beta({\gamma}_{pm+i})}{x}-\binom{\# \beta({\gamma}_{pm+i})
}{2}-(\sum \limits_{y\in \beta({\gamma}_{pm+i-1})}{y}-\binom{\#
\beta({\gamma}_{pm+i-1}) }{2})
 \\ &=&
  \sum
\limits_{1\leq k \leq m+1} (k t-p(m+1-k)-i)- \sum \limits_{1\leq k
\leq m+1}(\# \beta({\gamma}_{pm+i})-k)\\ &=&
\frac{m+1}{2}((m+2)t-pm-2i)- \sum \limits_{1\leq k \leq
m+1}(p\binom{m+1}{2}+(m+1)i-k) \\ &=& \frac{m+1}{2}((m+2)t-pm-2i)-
\frac{m+1}{2} (2p\binom{m+1}{2}+2(m+1)i-m-2)\\ &=& \frac{m+1}{2}
                 (m+2)(t-pm-2i+1) \\ &=& \binom{m+2}{2}(t-pm-2i+1).
\end{eqnarray*}

\textbf{(2)} By (1) we know
\begin{eqnarray*} f(pm+i)-f(pm) &=&
\sum \limits_{1\leq l \leq i}(f(pm+l)-f(pm+l-1))
\\
&=& \sum \limits_{1\leq l \leq i} \binom{m+2}{2}(t-pm-2l+1) \\
&=&
 \binom{m+2}{2}(it-ipm-i^2). \    \end{eqnarray*}
\hfill  $\square$

\noindent\textbf{Remark.} Notice that Lemma
\ref{f(pm+i+1)-f(pm+i)}(2) is also true for $i=0$.

\begin{lemma} \label{pm+i}
Let $0\leq i \leq p$ and $m\geq 0.$ Suppose that $pm+i\leq t-1.$
Then  $$f(pm+i)=\binom{m+2}{2}(it-i p m-i^2)
+\binom{m+2}{3}(pt-p^2)-3\binom{m+2}{4}p^2.$$
\end{lemma}
\textbf{Proof.}      By Lemma \ref{f(pm+i+1)-f(pm+i)} we have
\begin{eqnarray*}
f(pm+i)&=& f(pm+i)-f(pm)+\sum \limits_{0\leq k\leq m-1}(f(pk+p)-f(p
k))\\ &=& \binom{m+2}{2}(it-i p m-i^2) +\sum \limits_{0\leq k\leq
m-1}\binom{k+2}{2}(pt-p^2k-p^2) \\  &=& \binom{m+2}{2}(it-i p m-i^2)
+\binom{m+2}{3}(pt-p^2)-3\binom{m+2}{4}p^2.
 \end{eqnarray*}
In the above proof, we use the identities
$$k\binom{k+2}{2}=3\binom{k+2}{3}$$ and $$\sum_{0\leq k\leq m-1}
\binom{k+n}{n}=\binom{n+m}{n+1}.$$ \hfill $\square$

%By Lemma \ref{a-mt} and Lemma \ref{a_k}, there exists a unique $(t,
%t + 1, t + 2)$-core partition ${\lambda}_0$ such that
%$\beta({\lambda}_0)= \bigcup_{1\leq k \leq [\frac{t}{2}]} S_{t,k}$.

%\begin{lemma} \label{0}
%The size of ${\lambda}_0$ is $\mid{\lambda}_0\mid=   m
%\binom{m+1}{3}$ if $t = 2m-1$; $( m+1)
%\binom{m+1}{3}+\binom{m+2}{3}$ if $t = 2m$.
%\end{lemma}
%\textbf{Proof.}   Direct computation.

%%S_{t,k}}x=\sum \limits_{x=(k-1)(t+2)+1}^{kt-1}=\frac{(t-2k+1)((2k-1)t+2k-2)}{2}=
%\frac{(2k-1)t^2-(4k^2-6k+3)t-4k^2+6k-2}{2}$.

%Then we have $\mid{\lambda}_0\mid= \sum \limits_{1\leq k \leq
%[\frac{t}{2}]}\sum \limits_{x\in S_{t,k}}x-\binom{\sum \limits_{1\leq k \leq
%[\frac{t}{2}]}a_{t,k}}{2}= \sum \limits_{1\leq k \leq [\frac{t}{2}]}
%\frac{(2k-1)t^2-(4k^2-6k+3)t-4k^2+6k-2}{2}-\binom{\sum \limits_{1\leq k \leq
%[\frac{t}{2}]}(t-2k+1)}{2}=$

Now we prove our main result in this paper.

\vspace{2ex}

\noindent\textbf{Proof of Theorem \ref{main}.} Suppose that
$\lambda$ is a $(t,t+1,\ldots, t+p)$-core partition with the largest
size. By Lemma \ref{set}, we know $\beta(\lambda)\subseteq
\bigcup_{1\leq k \leq [\frac{t+p-2}{p}]} S_{k}$.  We will give some
properties for such $\lambda$.

\textbf{Step 1.} Let $r$ be the largest positive integer $k$ such
that $b_{\lambda,k}> 0$, i.e.,  $b_{\lambda,r}> 0$ and
$b_{\lambda,k}=0 $ for $ k> r$.  We claim that $\lambda=\mu_{c_1,
c_2, \ldots, c_r}$ for some positive integers $c_1, c_2, \ldots,
c_r$ such that $c_k\leq a_k$ for $1\leq k \leq r $ and  $
c_k-c_{k+1}\geq p$ for $1\leq k \leq r-1 $:

First we know $b_{\lambda,k}-b_{\lambda,k+1}\geq p$ for $1\leq k
\leq r-1 $ by Lemma \ref{bk}. It follows that $$kt-b_{\lambda,k}\leq
(k+1)t-b_{\lambda,k+1}-(t+p).$$ This means that for every
$$(k+1)t-b_{\lambda,k+1}\leq x \leq  (k+1)t-1 ,$$ we have
$$kt-b_{\lambda,k}\leq x-(t+i) \leq  kt-1 $$ for $0\leq i \leq p.$
Thus by Lemma \ref{a-mt} we have $ \bigcup_{1\leq k \leq r} \{ x\in
\textbf{Z} : kt-b_{\lambda,k}\leq x \leq kt-1 \}$ must be a
$\beta$-set for some $(t,t+1,\ldots, t+p)$-core partition
${\lambda}'$. We can write $$\beta({\lambda}')= \bigcup_{1\leq k
\leq r} \{ x\in \textbf{Z} : kt-b_{\lambda,k}\leq x \leq kt-1 \}.$$
Since $\{ x\in \textbf{Z} : kt-b_{\lambda,k}\leq x \leq  kt-1 \}$ is
just the set of the largest $b_{\lambda,k}$ elements in $S_{k}$, we
have
$$\sum \limits_{x\in \beta(\lambda)\bigcap S_{k}}x\ \leq \sum
\limits_{kt-b_{\lambda,k}\leq x \leq kt-1}x\ =\sum \limits_{x\in
\beta({\lambda}')\bigcap S_{k}}x$$ for $1\leq k\leq r.$ Then we know
$$ \mid \lambda \mid= \sum \limits_{x\in
\beta(\lambda)}x-\binom{\sum \limits_{1\leq k \leq r}
b_{\lambda,k}}{2} \leq  \sum \limits_{x\in
\beta({\lambda}')}x-\binom{\sum \limits_{1\leq k \leq r}
b_{\lambda,k}}{2} = \mid \lambda' \mid.$$ The above equality holds
if and only if $\lambda= \lambda'$. Since we already assume that
$\lambda$ is a $(t,t+1,\ldots, t+p)$-core partition with the largest
size, we must have $\lambda= \lambda'$ and thus
$$\beta(\lambda)\bigcap S_{k}=\{ x\in \textbf{Z} :
kt-b_{\lambda,k}\leq x \leq  kt-1 \}$$ for $1\leq k \leq r$. Let
$c_k=b_{\lambda,k}$ for $1\leq k\leq r$. Then we have
$\lambda=\mu_{c_1, c_2, \ldots, c_r}$, $c_k\leq a_k$ for $1\leq k
\leq r $ and $ c_k-c_{k+1}\geq p$ for $1\leq k \leq r-1$.
 We prove this
claim.

\textbf{Step 2.}  We claim that $1\leq c_r\leq p$:

Otherwise, suppose that $c_r\geq p+1$. Then $$r\leq
[\frac{t+p-2}{p}]-1$$ since $$c_{[\frac{t+p-2}{p}]}\leq
a_{[\frac{t+p-2}{p}]}\leq p.$$  Then we know $1 \leq a_{r+1}.$ Thus
we can define
$$\lambda'=\mu_{c_1, c_2, \ldots, c_r,  1}.$$ By Lemma
\ref{mostnumbers}, $\lambda'$ is a $(t,t+1,\ldots, t+p)$-core
partition since $c_r-1\geq p$. It is easy too see that
$$\beta({\lambda}')= \beta({\lambda})\bigcup \{ (r+1)t-1\}.$$ But by
Lemma \ref{thm1} we have \begin{eqnarray*} \mid \lambda' \mid-\mid
\lambda \mid &=& \sum \limits_{x\in \beta(\lambda')}{x}-\binom{\#
\beta(\lambda') }{2}-(\sum \limits_{y\in
\beta(\lambda)}{y}-\binom{\# \beta(\lambda) }{2})\\&=& (r+1)t-1-\#
\beta(\lambda) >0\end{eqnarray*} since $(r+1)t-1$ is larger than any
element in $\beta({\lambda})$. This contradicts the assumption that
$\lambda$ is a $(t, t + 1, \ldots, t + p)$-core partition with the
largest size. Then we must have $1\leq c_r\leq p$.

\textbf{Step 3.}  We claim that there is at most one integer $i$
satisfying $1\leq i \leq r-1$ and $ c_i-c_{i+1}\neq p$:

Otherwise, suppose that $1\leq i<j \leq r-1$ such that $
c_i-c_{i+1}\neq p$ and $ c_j-c_{j+1}\neq p$.   It is easy to see
that $$c_{j+1}+1\leq c_j-p \leq a_j-p= a_{j+1}.$$ Then we can define
$$\lambda'=\mu_{c_1, c_2, \ldots, c_{i-1}, c_i-1, c_{i+1},
\ldots, c_{j}, c_{j+1}+1, c_{j+2}, \ldots, c_r}.$$ By Lemma
\ref{mostnumbers}, $\lambda'$ is a $(t,t+1,\ldots, t+p)$-core
partition since $ c_i-c_{i+1}\geq p+1$ and $ c_j-c_{j+1}\geq p+1$.
It is easy too see that $$\beta({\lambda}')= \beta({\lambda})\bigcup
\{ (j+1)t-c_{ j+1}-1\}\setminus \{ it-c_i \}.$$ But by Lemma
\ref{thm1} we have $$\mid \lambda' \mid-\mid \lambda \mid=
(j+1)t-c_{ j+1}-1-(it-c_i)\geq 2t+c_i-c_{ j+1}-1>0.$$ This
contradicts the assumption that $\lambda$ is a $(t, t + 1, \ldots, t
+ p)$-core partition with the largest size. Then we prove the claim.

\textbf{Step 4.} We claim that if such $i$ in Step $3$ exists, then
$ c_i-c_{ i+1}= p+1$:

Otherwise, suppose that $ c_i-c_{ i+1}\geq p+2$.     It is easy to
see that $$c_{i+1}+1< c_i-p \leq a_i-p= a_{i+1}.$$ Then we can
define
$$\lambda'=\mu_{c_1, c_2, \ldots, c_{i-1}, c_i-1,
c_{i+1}+1, c_{i+2},  \ldots, c_r}.$$ By Lemma \ref{mostnumbers},
$\lambda'$ is a $(t,t+1,\ldots, t+p)$-core partition since $$
(c_i-1)-(c_{i+1}+1)=c_i-c_{i+1}-2\geq p.$$ Notice that
$$\beta({\lambda}')= \beta({\lambda})\bigcup \{
(i+1)t-c_{i+1}-1\}\setminus \{ it-c_{i} \}.$$ By Lemma \ref{thm1} we
have $$\mid \lambda' \mid-\mid \lambda \mid=
(i+1)t-c_{i+1}-1-(it-c_{i})=t+c_i-c_{ i+1}-1>0.$$ This contradicts
the assumption that $\lambda$ is a $(t, t + 1, \ldots, t + p)$-core
partition with the largest size. Then we prove the claim.

\textbf{Step 5.} We claim that such $i$ in Step 4 couldn't exist,
i.e., there is no such $i$ satisfying $ c_i-c_{i+1}=p+1$:

Otherwise, suppose that $ c_i-c_{i+1}=p+1$. It is easy to see that
$$c_{i+1}+1= c_i-p \leq a_i-p= a_{i+1}.$$ Then we can define
$$\lambda'=\mu_{c_1, c_2, \ldots, c_{i}, c_{i+1}+1, c_{i+2},
\ldots, c_r}$$ and $$\lambda''=\mu_{c_1, c_2, \ldots, c_{i-1},
c_{i}-1, c_{i+1}, \ldots, c_r}.$$ Then $\lambda'$ and $\lambda''$
are also $(t, t + 1, \ldots, t + p)$-core partitions by $
c_i-c_{i+1}=p+1$ and Lemma \ref{mostnumbers}. Notice that
$$\beta({\lambda}')= \beta({\lambda})\bigcup \{ (i+1)t-c_{i+1}-1
\}$$ and $$\beta({\lambda}'')= \beta({\lambda})\setminus \{ it-c_{i}
\}.$$ By Lemma \ref{thm1} we have
$$\mid \lambda' \mid-\mid \lambda \mid=(i+1)t-c_{i+1}-1-\# \beta(\lambda)=(i+1)t-c_{i+1}-1-\sum \limits_{1\leq k \leq r}c_k$$ and

$$\mid \lambda \mid-\mid \lambda'' \mid
=it-c_{i}-\# \beta(\lambda'') =it-c_{i}-(\sum \limits_{1\leq k \leq
r}c_k-1).$$ Put these two equalities together, we have
\begin{eqnarray*}
2\mid \lambda \mid  & = & \mid \lambda'\mid+\mid
\lambda''\mid-(t+c_i-c_{i+1}-2)
\\ &=&
 \mid \lambda'\mid+\mid
\lambda''\mid-(t+p-1) \\ &<& \mid \lambda'\mid+\mid \lambda''\mid.
\end{eqnarray*}

This contradicts the assumption that $\lambda$ is a $(t, t + 1,
\ldots, t + p)$-core partition with  the largest size. Then we prove
the claim.

\textbf{Step 6.} We claim that $\lambda\in \{ {\gamma}_{t-j}:\ 1\leq
j \leq p \}$:

By Step $2,\ 3,\ 4,\ 5$ we know $1\leq c_r \leq p$ and  $
c_k-c_{k+1}= p$ for $1\leq k \leq r-1,$ which means that $$\lambda=
\mu_{c_1, c_1-p, c_1-2p,  \ldots, c_1- [\frac{c_1-1}{p}]
p}={\gamma}_{c_1}.$$ Suppose that $c_1=pm+i$, where $0\leq i\leq
p-1$ and $m\geq 0.$ If $c_1<t-p$, then $\gamma_{c_1+1}$ is well
defined and by Lemma \ref{f(pm+i+1)-f(pm+i)} we have
\begin{eqnarray*} \mid \gamma_{c_1+1} \mid -\mid \lambda \mid & = & \mid \gamma_{c_1+1} \mid -\mid \gamma_{c_1}\mid = f(c_1+1)-f(c_1)\\
&=& f(pm+i+1)-f(pm+i) = \binom{m+2}{2}(t-pm-2(i+1)+1)\\&=&
\binom{m+2}{2}(t-c_1-i-1) > \binom{m+2}{2}(p-i-1)\geq 0.
\end{eqnarray*} This means that  $c_1<t-p$ implies $\mid
\gamma_{c_1+1} \mid > \mid \lambda \mid$, which contradicts the
assumption that $\lambda$ is a $(t, t + 1, \ldots, t + p)$-core
partition with the largest size. Then we must have $c_1\geq t-p$.
But $c_1\leq a_1=t-1$, thus $\lambda\in \{ {\gamma}_{t-j}:\ 1\leq j
\leq p \}$.

\textbf{Step 7.} By assumption we know $t=pn+d$, where $1\leq d\leq
p$ and $n\geq 0.$  We claim that $\mid \lambda \mid =
\text{max}\{f(pn+[\frac{d}{2}]),\ f(pn-[\frac{p-d}{2}]) \}$ and
$\lambda\in \{ {\gamma}_{j}:\ j=pn+[\frac{d}{2}],\
pn+[\frac{d}{2}]+1,\ pn-[\frac{p-d}{2}],\ \text{or} \
pn-[\frac{p-d}{2}]-1\}$:

By Step $6$ we know $$\lambda\in \{ {\gamma}_{pn+k}:\ 0\leq k \leq
d-1\}\bigcup \{ {\gamma}_{pn-k}:\ 0\leq k \leq p-d\}.$$

For ${\gamma}_{pn+k}$ where $0\leq k \leq d-1$, by Lemma
\ref{f(pm+i+1)-f(pm+i)} we have
\begin{eqnarray*}f(pn+k)-f(pn) &=&\binom{n+2}{2}(kt-kpn-k^2)\\ &=&\binom{n+2}{2}(kd-k^2)\\ &=&\binom{n+2}{2}(\frac{d^2}{4}-(k-\frac{d}{2})^2).\end{eqnarray*}
Then it is easy to see that when $d$ is even, $f(pn+k)$ is maximal
for $0\leq k \leq d-1$ if and only if $k=[\frac{d}{2}]$; when $d$ is
odd, $f(pn+k)$ is maximal for $0\leq k \leq d-1$ if and only if
$k=[\frac{d}{2}]$ or $ [\frac{d}{2}]+1.$

For ${\gamma}_{pn-k}$ where $ 0\leq k \leq p-d$,  by Lemma
\ref{f(pm+i+1)-f(pm+i)} we have
\begin{eqnarray*}f(pn)-f(pn-k)&=&f(p(n-1)+p)-f(p(n-1)+p-k)\\ &=&\sum \limits_{p-k\leq l \leq p-1
}(f(p(n-1)+l+1)-f(p(n-1)+l) )\\ &=&\sum \limits_{p-k\leq l \leq p-1
}\binom{n+1}{2}(t-p(n-1)-2(l+1)+1)\\ &=& \sum \limits_{p-k\leq l
\leq
p-1 }\binom{n+1}{2}(p+d-2l-1)\\
&=&\binom{n+1}{2}(k^2-k(p-d))\\
&=&\binom{n+1}{2}((k-\frac{p-d}{2})^2-\frac{(p-d)^2}{4}).\end{eqnarray*}
Then it is easy to see that when $p-d$ is even, $f(pn-k)$ is maximal
for $ 0\leq k \leq p-d$ if and only if $k=[\frac{p-d}{2}]$; when
$p-d$ is odd, $f(pn-k)$ is maximal for $ 0\leq k \leq p-d$ if and
only if $k=[\frac{p-d}{2}]$ or $ [\frac{p-d}{2}]+1.$ Then we prove
the claim.

\textbf{Step 8.}  By Step $7$ and Lemma \ref{pm+i} we have the
largest size of $(t, t + 1, \ldots, t + p)$-core partitions is
$$\text{max}\{f(pn+[\frac{d}{2}]),\ f(pn-[\frac{p-d}{2}]) \},$$
which equals to
\begin{eqnarray*}
&\max\{\binom{n+2}{2}[\frac{d}{2}](d-[\frac{d}{2}])+\binom{n+2}{3}(p^2n+pd-p^2)-3\binom{n+2}{4}p^2,&\\
&\binom{n+1}{2}(p-[\frac{p-d}{2}])(d+[\frac{p-d}{2}])+\binom{n+1}{3}(p^2n+pd-p^2)-3\binom{n+1}{4}p^2
\}.&\end{eqnarray*}

By Step $7$ we also know the number of $(t, t + 1, \ldots, t +
p)$-core partitions with the largest size is at most $4$. We finish
the proof. \hfill $\square$

For $p=1$, we have the following corollary. We mention that, this
corollary could be obtained by results on the  largest size of
$(t_1,  t_2)$-core partitions  in \cite{ols}.

\begin{cor} \label{p=1}
The  largest size of $(t, t + 1)$-core partition is $
\binom{t+2}{4}$. The number of $(t, t + 1)$-core partitions with the
largest size is $1$.
\end{cor}
\textbf{Proof.} By Step $6$ of Theorem \ref{main}, we know
${\gamma}_{t-1}$ is the only $(t, t + 1)$-core partition with the
largest size since $p=1$ in this case.

By Lemma \ref{pm+i} the largest size is

$$
\mid {\gamma}_{t-1} \mid  = f(t-1) = \binom{t+1}{3}(t-1)
-3\binom{t+1}{4}= \binom{t+2}{4}.$$  \hfill $\square$

For $p=2$, we have the following corollary, which shows that
Amdeberhan's conjecture on the largest size of $(t,t+1, t+2)$-core
partitions proposed in \cite{tamd} is true.

\begin{cor} \label{p=2}
\textbf{(1)} If $\ t=2n-1$, the  largest size of $(2n-1, 2n,
2n+1)$-core partition is $n \binom{n+1}{3}.$ The number of $(2n-1,
2n, 2n+1)$-core partitions with the largest size is $2$.

\textbf{(2)} If $\ t=2n$, the  largest size of $(2n, 2n + 1, 2n +
2)$-core partition is $( n+1) \binom{n+1}{3}+\binom{n+2}{3}$. The
number of $(2n, 2n + 1, 2n + 2)$-core partitions with the largest
size is $1$.
\end{cor}
\textbf{Proof.} By Step $6$ of Theorem \ref{main}, we know a $(t, t
+ 1,  t + 2)$-core partition with the largest size must be
${\gamma}_{t-1}$ or ${\gamma}_{t-2}$ since $p=2$ in this case.

\textbf{(1)} When $t=2n-1$, by Lemma \ref{f(pm+i+1)-f(pm+i)} we have
$$\mid {\gamma}_{t-1} \mid- \mid {\gamma}_{t-2} \mid=
f(2n-2)-f(2n-3)=\binom{n}{2}( t-2(n-2)-4+1 )= 0.$$ Then we get
$f(2n-2)=f(2n-3)$. This means that ${\gamma}_{2n-2}$ and
${\gamma}_{2n-3}$ are the only two $(2n-1, 2n, 2n+1)$-core
partitions with the largest size. By Lemma \ref{pm+i} the largest
size is
\begin{eqnarray*}
f(2n-2) &=& f(2(n-1))\\ &=&\binom{n+1}{3}(2(2n-1)-4)-12\binom{n+1}{4}\\
&=& \binom{n+1}{3}(4n-6-3(n-2))  \\ &=& n\binom{n+1}{3}.
\end{eqnarray*}

\textbf{(2)} When $t=2n$, by Lemma \ref{f(pm+i+1)-f(pm+i)} we have
\begin{eqnarray*}\mid {\gamma}_{t-1} \mid- \mid {\gamma}_{t-2}
\mid &=& f(2n-1)-f(2n-2)\\ &=& \binom{n+1}{2}( t-2(n-1)-2+1
)\\&=&\binom{n+1}{2}> 0.\end{eqnarray*} This means that
${\gamma}_{2n-1}$ is the only $(2n, 2n + 1, 2n + 2)$-core partition
with  the largest size. By Lemma \ref{pm+i} the largest size is
\begin{eqnarray*}
f(2n-1) &=& f(2(n-1)+1)\\ &=& \binom{n+1}{2}(2n- 2 (n-1)-1)
+\binom{n+1}{3}(4n-4)-12\binom{n+1}{4}\\ &=& \binom{n+1}{2}
+\binom{n+1}{3}(4n-4)-3\binom{n+1}{3}(n-2)\\ &=& \binom{n+1}{2}
+(n+2)\binom{n+1}{3}\\ &=& \binom{n+1}{2}
+\binom{n+1}{3}+(n+1)\binom{n+1}{3}\\
&=& (n+1)\binom{n+1}{3}+\binom{n+2}{3}.
\end{eqnarray*}
 \hfill $\square$

\section{Acknowledgements}
We'd like to thank Prof. P. O. Dehaye for the helpful conversations
and discussions. The author is supported  by Forschungskredit of the
University of Zurich, grant no. [FK-14-093].


\begin{thebibliography}{1}

\bibitem{tamd}
T. Amdeberhan, Theorems, problems and conjectures, $2013.$ Published
electronically at
www.math.tulane.edu/$\sim$tamdeberhan/conjectures.pdf.

\bibitem{and}
J. Anderson, Partitions which are simultaneously $t_1$- and
$t_2$-core, Disc. Math. $248(2002),\ 237-243$.


\bibitem{berge} C. Berge, Principles of Combinatorics, Mathematics in Science and Engineering Vol. $72$, Academic
Press, New York, $1971$.

\bibitem{Herman}
 F. Chung and J. Herman, Some results on hook lengths, Discrete Math.
$20  (1977), \ 33-40$.



\bibitem{stanton2}
F. Garvan, D. Kim, and D. Stanton, Cranks and $t$-cores, Inv. Math.
$101 (1990),\ 1-17.$

\bibitem{gran}
A. Granville and K. Ono, Defect zero $p$-blocks for finite simple
groups, Trans. Amer. Math. Soc. $348 (1996),\ 331-347.$


\bibitem{xiong}
H. Xiong, The number of simultaneous core partitions,
arXiv:$1409.7038$[math.CO], $2014$.


\bibitem{jk}
G. James and A. Kerber, The Representation Theory of the Symmetric
Group, Addison-Wesley Publishing Company, Reading, MA, $1981$.

%\bibitem{frame}
%J. Frame,  G. Robinson,  and R. Thrall, The hook graphs of the
%symmetric group, Canad. J. Math. $6 (1954), \ 316-324.$


\bibitem{ols}
J. Olsson and D. Stanton, Block inclusions and cores of partitions,
Aequationes Math. $74(1-2)(2007),\ 90-110$.

\bibitem{stanton}
D. Stanton, Open positivity conjectures for integer partitions,
Trends Math. $2 (1999),\ 19 - 25$.






\end{thebibliography}
\end{document}